\documentclass[12pt]{article}
\usepackage{amsmath,amssymb,amsthm,amsfonts,amsbsy, enumerate}
\usepackage[dvips]{epsfig}

\theoremstyle{plain}
\newtheorem{theorem}{Theorem}[section]
\newtheorem{corollary}[theorem]{Corollary}
\newtheorem{conj}[theorem]{Conjecture}
\newtheorem{lemma}[theorem]{Lemma}

\newcommand{\mc}{{\rm mc}}

\begin{document}

\title{Max-cut and extendability of matchings in distance-regular graphs}
\author{Sebastian M. Cioab\u{a}\thanks{Research supported in part by NSA grant H98230-13-1-0267 and NSF grant DMS-1600768.}\\
\small Department of Mathematical Sciences\\[-0.8ex]
\small University of Delaware\\[-0.8ex] 
\small Newark, DE 19707-2553, U.S.A.\\
\small\tt cioaba@udel.edu\\
\and
Jack Koolen\thanks{Research partially supported by the 100 talents program of Chinese Academy of Sciences. Research is also partially supported by the National Natural Science Foundation of China (No. 11471009). }\\
\small  Wu Wen-Tsun Key Laboratory of Mathematics of CAS \\[-0.8ex]
\small  School of Mathematical Sciences \\[-0.8ex]
\small  University of Science and Technology of China \\[-0.8ex]
\small  96 Jinzhai Road, Hefei, 230026, Anhui, P.R. China \\
\small\tt koolen@ustc.edu.cn\\
\and
Weiqiang Li\thanks{Research supported by the University Dissertation Fellows award by University of Delaware.}  \\
\small Department of Mathematical Sciences\\[-0.8ex]
\small University of Delaware\\[-0.8ex]
\small Newark, DE 19707-2553, U.S.A.\\
\small\tt weiqiang@udel.edu
}
\date{\today}

\maketitle

\begin{abstract}
A connected graph $G$ of even order $v$ is called $t$-extendable if it contains a perfect matching, $t<v/2$ and any matching of $t$ edges is contained in some perfect matching. The extendability of $G$ is the maximum $t$ such that $G$ is $t$-extendable. Since its introduction by Plummer in the 1980s, this combinatorial parameter has been studied for many classes of interesting graphs. In 2005, Brouwer and Haemers proved that every distance-regular graph of even order is $1$-extendable and in 2014, Cioab\u{a} and Li showed that any connected strongly regular graph of even order is $3$-extendable except for a small number of exceptions.

In this paper, we extend and generalize these results. We prove that all distance-regular graphs with diameter $D\geq 3$ are $2$-extendable and we also obtain several better lower bounds for the extendability of distance-regular graphs of valency $k\geq 3$ that depend on $k$, $\lambda$ and $\mu$, where $\lambda$ is the number of common neighbors of any two adjacent vertices and $\mu$ is the number of common neighbors of any two vertices in distance two. In many situations, we show that the extendability of a distance-regular graph with valency $k$ grows linearly in $k$. We conjecture that the extendability of a distance-regular graph of even order and valency $k$ is at least $\lceil k/2\rceil-1$ and we prove this fact for bipartite distance-regular graphs.

In course of this investigation, we obtain some new bounds for the max-cut and the independence number of distance-regular graphs in terms of their size and odd girth and we prove that our inequalities are incomparable with known eigenvalue bounds for these combinatorial parameters.
\end{abstract}

\section{Introduction}

Our graph theoretic notation is standard (for undefined notions, see \cite{BH2, GR,VW}). The adjacency matrix of a graph $G=(V,E)$ has its rows and columns indexed after the vertices of the graph and its $(u,v)$-th entry equals $1$ if $u$ and $v$ are adjacent and 0 otherwise. If $G$ is a connected $k$-regular graph of order $v$, then $k$ is the largest eigenvalue of the adjacency matrix of $G$ and its multiplicity is 1. In this case, let $k=\lambda_1>\lambda_2\geq \dots \geq \lambda_v$ denote the eigenvalues of the adjacency matrix of  $G$. If $S$ and $T$ are vertex disjoint subsets of $G$, let $e(S,T)$ denote the number of edges with one endpoint in $S$ and the other in $T$. If $S$ is a subset of vertices of $G$, let $S^c$ denote its complement. The max-cut of $G$ is defined as $\mc(G):=\max_{S\subset V} e(S,S^c)$ and measures how close is $G$ from being a bipartite graph. Given a graph $G$, determining $\mc(G)$ is a well-known NP-hard problem (see \cite[Problem ND16, page 210]{GJ} or \cite{Karp}) and designing efficient algorithms to approximate $\mc(G)$ has attracted a lot of attention \cite{AS,DP1,DP2,DP3,GR,GW,MP,Tre}.
 
A set of edges $M$ in a graph $G$ is a matching if no two edges of $M$ share a vertex. A matching $M$ is perfect if every vertex is incident with exactly one edge of $M$. A connected graph $G$ of even order $v$ is called $t$-extendable if it contains at least one perfect matching, $t<v/2$ and any matching of size $t$ is contained in some perfect matching. Graphs that are $1$-extendable are also called matching-covered (see Lov\'{a}sz and Plummer \cite[page 113]{LP}). The extendability of a graph $G$ of even order is defined as the maximum $t<v/2$ such that $G$ is $t$-extendable. This concept was introduced by Plummer \cite{Plummer80} in 1980 and was motivated by work of Lov\'{a}sz \cite{L1} on canonical decomposition of graphs containing perfect matchings. Later on, Yu \cite{Yu} expanded the definition of extendability to graphs of odd order. Zhang and Zhang \cite{ZZ} obtained an $O(mn)$ algorithm to compute the extendability of a bipartite graph with $n$ vertices and $m$ edges, but the complexity of determining the extendability of a general graph is unknown at present time (see \cite{P1,P2,YuLiu} for more details on extendability of graphs).

In this paper, we obtain a simple upper bound for the max-cut of certain regular graphs in terms of their odd girth (the shortest length of an odd cycle). In Section 2, we prove that if $G$ is a non-bipartite distance-regular graph with $e$ edges and odd girth $g$, then $\mc(G) \leq e(1-\frac{1}{g})$. As a consequence of this result, we show that if $G$ is a non-bipartite distance-regular graph with $v$ vertices, odd girth $g$ and independence number $\alpha(G)$, then $\alpha(G) \leq \frac{v}{2} (1-\frac{1}{g})$. We show that these bounds are incomparable with some spectral bounds of Mohar and Poljak \cite{MP} for the max-cut and of Cvetkovi\'{c} (see \cite[Theorem 3.5.1]{BH2} or \cite[Lemma 9.6.3]{GR}) and Hoffman (see \cite[Theorem 3.5.2]{BH2} or \cite[Lemma 9.6.2]{GR}) for the independence number.

Holton and Lou \cite{HL} showed that strongly regular graphs with certain connectivity properties are $2$-extendable and conjectured that all but a few strongly regular graphs are $2$-extendable. Lou and Zhu \cite{LZ} proved this conjecture and showed that every connected strongly regular graph of valency $k\geq 3$ is $2$-extendable with the exception of the complete $3$-partite graph $K_{2,2,2}$ and the Petersen graph. Cioab\u{a} and Li \cite{CL} proved that every connected strongly regular graph of valency $k\geq 5$ is $3$-extendable with the exception of the complete $4$-partite graph $K_{2,2,2,2}$, the complement of the Petersen graph and the Shrikhande graph. Moreover, Cioab\u{a} and Li determined the extendability of many families of strongly regular graphs including Latin square graphs, block graphs of Steiner systems, triangular graphs, lattice graphs and all known triangle-free strongly regular graphs. 
For any such graph of valency $k$, Cioab\u{a} and Li proved that the extendability is at least $\lceil k/2 \rceil -1$ and conjectured that this fact should be true for any strongly regular graph. 

In this paper, we extend and generalize these results and study the extendability of distance-regular graphs with diameter $D\geq 3$. Brouwer and Haemers \cite{BH} proved that distance-regular graphs are $k$-edge-connected. Plesn\'{i}k (\cite{Plesnik} or \cite[Chapter 7]{L}) showed that if $G$ is a $k$-regular $(k-1)$-edge-connected graph with an even number of vertices, then the graph obtained by removing any $k-1$ edges of $G$ contains a perfect matching. These facts imply that every distance-regular graph of even order is $1$-extendable. In Section 3, we improve this result and we show that all distance-regular graphs with diameter $D\geq 3$ are $2$-extendable. We prove that any distance-regular graph of valency $k\geq 3$ with $\lambda \geq 1$ is $\lfloor \frac{k+1-\frac{k}{\lambda+1}}{2}\rfloor $-extendable (when $\mu=1$), $\lfloor \frac{1}{2}\lceil \frac{k+2}{2} \rceil \rfloor$-extendable (when $\mu=2$) and $\lfloor \frac{k}{3} \rfloor$-extendable (when $\mu \geq 3$ and $k\geq 6$). We also show that any bipartite distance-regular graph of valency $k$ is $\lfloor \frac{k+1}{2} \rfloor$-extendable. We also remark that our results for graphs of even order can be extended to graphs of odd order in similar fashion to what was done in \cite{CL}, but for the sake of brevity and clarity, we will not include the details here.

\section{Max-cut of distance regular graphs}

For notation and definitions related to distance-regular graphs,
see \cite{BCN}. We denote the intersection array of a distance-regular $G$ of diameter $D$ by $\{b_0,\ldots,b_{D-1};\,c_1,\ldots,c_D\}$ and we let $k = b_0$ and $a_i=k-b_i-c_i$ for $0\leq i\leq D$ as usual. Also, let $\lambda=a_1$ and $\mu=c_2$. The following result gives a simple upper bound for the max-cut of a graph in terms of its odd girth under certain regularity conditions. Such regularity conditions will be satisfied by walk-regular graphs and distance-regular graphs. As pointed to us by one of the anonymous referees, the theorem below holds for any odd natural number $g$ as long as the condition that every edge is in the same number of cycles of length $g$, is satisfied.
\begin{theorem}\label{maxcut}
Let $G$ be a non-bipartite graph with odd girth $g$. If every edge of $G$ is contained in the same number of cycles of length $g$, then 
 \begin{equation}\label{maxcutbound}		
 \mc(G) \leq e\left (1-\frac{1}{g}\right ).
 \end{equation}
 \end{theorem}
\begin{proof}
Let $\gamma$ be the number of cycles of length $g$ containing some fixed edge of $G$ and let $\mathcal{C}$ be the set of cycles of length $g$. By counting pairs $(e_0,C)$ with $e_0\in E(G), \ C\in \mathcal{C}$ with $e_0$ contained in $C$, we get that $|\mathcal{C}|=\frac{e\gamma}{g}$. Let $A$ be any subset of vertices and $T$ be the set of the edges with both endpoints in $A$ or in $A^c$. Every time we delete an edge in $T$, we destroy at most $\gamma$ cycles in $\mathcal{C}$. Therefore $|T|\geq \frac{|\mathcal{C}|}{\gamma}=\frac{e}{g}$. Since $e(A,A^c)=e-|T| \leq e(1-\frac{1}{g})$, this implies the desired conclusion.
\end{proof}
Our theorem can be applied to the family of $m$-walk regular graphs with $m\geq 1$. This family of graphs contains the distance-regular graphs. A connected graph $G$ is $m$-walk-regular if the number of walks of length $l$ between any pair of vertices only depends on the distance between them, provided that this distance does not exceed $m$. The family of $m$-walk-regular graphs was first introduced by Dalf\'{o}, Fiol, and Garriga \cite{DFG, FG}.

Note that the upper bound of Theorem \ref{maxcut} is tight as shown for example by the blow up of an odd cycle $C_g$. Such a graph can be constructed from the odd cycle $C_g$ by replacing each vertex $i$ of $C_g$ by a coclique $A_i$ of size $m$ for $1\leq i \leq g$ and adding all the possible edges between $A_i$ and $A_j$ whenever $i$ and $j$ are adjacent in $C_g$. The resulting graph which is also the lexicographic product of the cycle $C_g$ with the empty graph of order $m$ (see \cite[Ex 26, p.17]{GR} for a definition), has $gm$  vertices and $gm^2$ edges. The odd girth of this graph is $g$, each edge of the graph is contained in the same number of cycles of length $g$ and there is a cut of size $e\left (1-\frac{1}{g}\right)=(g-1)m^2$.  

Mohar and Poljak \cite{MP} showed that $mc(G)\leq \frac{v\mu_{max}}{4}$ for any graph $G$ on $v$ vertices whose largest Laplacian eigenvalue is $\mu_{max}$ (see also \cite{AS,DS,DP1,DP2,DP3} for related results). Translated to regular graphs, their result implies the following inequality:
\begin{equation}\label{ASmaxcutbound}
\mc(G)\leq \frac{e}{2}\left(1-\frac{\lambda_v}{k}\right).
\end{equation}

Note that the inequalities \eqref{maxcutbound} and \eqref{ASmaxcutbound} are incomparable. This fact can be seen by considering the complete graphs and the odd cycles, but we give other examples of distance-regular graphs below. Also, a simple calculation yields that inequality \eqref{maxcutbound} is better for graphs that in a spectral sense {\em are closer to being bipartite} (when $\lambda_v\leq -k(1-2/g)$ more precisely).
 
The Hamming graph $H(D,q)$ is the graph whose vertices are all the words of length $D$ over an alphabet of size $q$ with two words being adjacent if and only their Hamming distance is $1$. The graph $H(D,q)$ is distance-regular of diameter $D$,  has eigenvalues $(q-1)D-qi$  for $0\leq i\leq D$  and is bipartite when $q=2$ \cite[page 174]{BH2}. When $q\geq 3$, inequality \eqref{maxcutbound} always gives an upper bound $\frac{2e}{3}$. The upper bound from inequality \eqref{ASmaxcutbound} is $\frac{e}{2}(1+\frac{1}{q-1})$.  When $q=3$, \eqref{maxcutbound} is better. When $q\geq 5$, inequality \eqref{ASmaxcutbound} is better. When $q=4$, both inequalities give the same upper bound. 

The Johnson graph $J(n,m)$ is the graph whose vertices are the $m$-subsets of a set of size $n$ with two $m$-subsets being adjacent if and only if they have $m-1$ elements in common. The graph $J(n,m)$ is distance-regular with diameter $D=\min(m,n-m)$, eigenvalues $(m-i)(n-m-i)-i$, where $0\leq i\leq D$ \cite[page 175]{BH2}. Inequality \eqref{maxcutbound} always gives an upper bound $\frac{2e}{3}$. Inequality \eqref{ASmaxcutbound} is $\frac{e}{2}(1+\frac{D}{m(n-m)})$. When $\max (m, n-m) \geq 4$, inequality \eqref{ASmaxcutbound} is better and in the other cases ($m\in\{2, 3\}$ or $n-m\in \{1,2,3\}$), \eqref{maxcutbound} is better. 

In the following examples, we compare \eqref{maxcutbound} and \eqref{ASmaxcutbound} for other distance-regular graph with larger odd girth. 

\begin{enumerate}

\item The Dodecahedron graph \cite[page 417]{BCN} is a $3$-regular graph of order $20$ and size $30$. It has $\lambda_v=-\sqrt{5}$ and $g=5$. Inequality \eqref{maxcutbound} gives $\mc(G)\leq 24$ and inequality \eqref{ASmaxcutbound} gives $\mc(G) \leq 26$. 

\item The Coxeter graph \cite[page 419]{BCN} is a $3$-regular graph of order $28$ and size $42$. It has $\lambda_v=-\sqrt{2}-1\approx -2.414$ and $g=7$. Inequality \eqref{maxcutbound} gives $\mc(G)\leq 36$ and inequality \eqref{ASmaxcutbound} gives $\mc(G) \leq 37$. 

\item The Biggs-Smith graph \cite[page 414]{BCN} is a $3$-regular graph of order $102$ and size $153$. It has $\lambda_v\approx -2.532$ and $g=9$. Inequality \eqref{maxcutbound} gives $\mc(G)\leq 136$ and inequality \eqref{ASmaxcutbound} gives $\mc(G) \leq 141$. 

\item The Wells graph \cite[page 421]{BCN} is a $5$-regular graph of order $32$ and size $80$. It has $\lambda_v=-3$ and $g=5$. Inequality \eqref{maxcutbound} gives $\mc(G)\leq 64$ and inequality \eqref{ASmaxcutbound} gives $\mc(G) \leq 64$. 

\item The Hoffman-Singleton graph \cite[page 391]{BCN} is a $7$-regular graph of order $50$ and size $175$. It has $\lambda_v=-3$ and $g=5$. Inequality \eqref{maxcutbound} gives $\mc(G)\leq 140$ and inequality \eqref{ASmaxcutbound} gives $\mc(G) \leq 125$. 

\item The Ivanov-Ivanov-Faradjev graph \cite[page 414]{BCN} is a $7$-regular graph of order $990$ and size $3465$. It has $\lambda_v=-4$ and $g=5$. Inequality \eqref{maxcutbound} gives $\mc(G)\leq 2772$ and inequality \eqref{ASmaxcutbound} gives $\mc(G) \leq 2722$. 

\item The Odd graph $O_{m+1}$ \cite[page 259-260]{BCN} is the graph whose vertices are the $m$-subsets of a set with $2m+1$ elements, where two $m$-subsets are adjacent if and only if they are disjoint. Note that $O_3$ is Petersen graph. The graph $O_{m+1}$ is a distance-regular graph of valency $m+1$, order $v={2m+1 \choose m}$ and size $e=\frac{m+1}{2}{2m+1 \choose m}$. It has $\lambda_v=-m$ and $g=2m+1$. Inequality \eqref{maxcutbound} gives $\mc(G)\leq e(1-\frac{1}{2m+1})$ and inequality \eqref{ASmaxcutbound} gives $\mc(G) \leq e(1-\frac{1}{2m+2})$. 
\end{enumerate}

Theorem \ref{maxcut} can be used to obtain an upper bound for the independence number of certain regular graphs.  
\begin{corollary}\label{independencenumber}
Let $G$ be a non-bipartite regular graph with valency $k$ and odd girth $g$. If every edge of $G$ is contained in the same number of cycles of length $g$, then 
\begin{equation}\label{independencenumberbound}
\alpha(G) \leq \frac{v}{2}\left (1-\frac{1}{g}\right ).
\end{equation}
\end{corollary}
\begin{proof}
Let $S$ be an independent set of size $\alpha(G)$. Then $k\alpha(G)=e(S,S^c) \leq \frac{vk}{2} (1-\frac{1}{g})$ which implies the conclusion of the theorem.
\end{proof}

The Cvetkovi\'c inertia bound (see \cite[Theorem 3.5.1]{BH2} or \cite[Lemma 9.6.3]{GR}) states that if $G$ is a graph with $n$ vertices whose adjacency matrix has $n_{+}$ positive eigenvalues and $n_{-}$ negative eigenvalues, then
\begin{equation}\label{Cvetkovic}
\alpha(G)\leq \min(n-n_{-},n-n_{+}).
\end{equation}

The Hoffman-ratio bound (see \cite[Theorem 3.5.2]{BH2} or \cite[Lemma 9.6.2]{GR}) states that if $G$ is a $k$-regular graph with $v$ vertices, then
\begin{equation}\label{Hoffman}
\alpha(G)\leq \frac{v}{1+k/(-\lambda_v)}.
\end{equation} 

In the table below, we compare the bounds \eqref{independencenumberbound}, \eqref{Cvetkovic}  and \eqref{Hoffman} for some of the previous examples. When the bounds obtained are not integers, we round them below. The exact values of the independence numbers below were computed using Sage.

\begin{center}

\begin{tabular}{|c|cccc|}
\hline
Graph & $\alpha$ & \eqref{independencenumberbound} &  \eqref{Cvetkovic} & \eqref{Hoffman}\\
\hline
Dodecahedron & $8$ & $8$ & $8$ & $11$  \\
\hline
Coxeter & $12$ & $12$ & $13$ & $12$ \\
\hline
Biggs-Smith & $43$ & $45$ & $58$ & $46$ \\ 
\hline
Wells & $10$ & $12$ & $13$ & $12$ \\
\hline
Hoffman-Singleton & $15$ & $ 20$ & $21$ & $15$ \\
\hline
\end{tabular}

\end{center}

For the Hamming graph $H(D,q)$ with $D=2$ and $q\geq 3$, \eqref{independencenumberbound} is better than \eqref{Hoffman}. For the Hamming graph $H(D,q)$ with $D\geq 3$ and $q\geq 3$, \eqref{Hoffman} is better. For the Odd graph $O_{m+1}$, the inequalities \eqref{independencenumberbound} and \eqref{Hoffman} give the same bound that equals the independence number of $O_{m+1}$.

\section{Extendability of matchings in distance-regular graphs}

In this section, we will focus on the extendability of distance-regular graphs of even order. Similar results can be obtained for distance-regular graphs of odd order using the definition of extendability of Yu \cite{Yu}, but for the sake of simplicity we restrict ourselves to graphs of even order. A connected graph $G$ of odd order $v$ containing at least one matching of size $\frac{v-1}{2}$ (a near perfect matching)
is called $t$-near-extendable (or $t1/2$-extendable in the notation of Yu \cite{Yu}) if $t<\frac{v-1}{2}$ and for every vertex $x$, any matching of size $t$ that does not cover $x$, is contained in some near perfect matching that misses $x$. Graphs that are $0$-near-extendable are also called factor-critical or hypomatchable (see Lov\'{a}sz and Plummer \cite[page 89]{LP}).

In Subsection \ref{tools}, we describe the main tools which will be used in our proofs. In Subsection \ref{lowerbnds}, we give various lower bounds for the extendability of distance-regular graphs. In Subsection \ref{extend2}, we show that all  distance-regular graphs with diameter $D\geq 3$ are 2-extendable.

\subsection{Main tools}\label{tools}

Let $o(G)$ denote the number of components of odd order in a graph $G$. If $S$ is a subset of vertices of $G$, then $G-S$ denotes the subgraph of $G$ obtained by deleting the vertices in $S$. Let $N(T)$ denote the set of vertices outside $T$ that are adjacent to at least one vertex of $T$. When $T=\{x\}$, let $N(x)=N(\{x\})$. The distance $d(x,y)$ between two vertices $x$ and $y$ of a connected graph $G$ is the shortest length of a path between $x$ and $y$. If $x$ is a vertex of a distance-regular graph $G$, let $N_i(x)$ denote the set of vertices at distance $i$ from vertex $x$ and $k_i=|N_i(x)|$; the $i$-th subconstituent $\Gamma_i(x)$ of $x$ is the subgraph of $G$ induced by $N_i(x)$. 

\begin{theorem}[Brouwer and Haemers \cite{BH}]\label{kedge}
Let $G$ be a distance-regular graph of valency $k$. Then $G$ is $k$-edge-connected. Moreover, if $k>2$, 
then the only disconnecting sets of $k$ edges are the set of $k$ edges on a single vertex. 
\end{theorem}

\begin{theorem}[Brouwer and Koolen \cite{BK2}]\label{kconnected}
Let $G$ be a distance-regular graph of valency $k$. Then $G$ is $k$-connected. Moreover, if $k>2$, then the only disconnecting sets of $k$ vertices are the set of 
the neighbors of some vertex. 
\end{theorem}

\begin{lemma}\label{outedgeofsmallsets}
Let $G$ be a distance-regular graph with $k\geq 4$. If $A\subset V$ with $3\leq |A|\leq k-1$, then $e(A,A^c)\geq 3k-6$.
\end{lemma}
\begin{proof}
If $|A|\leq k-2$, then every vertex in $A$ has at least $k-(|A|-1)$ many neighbors in $A^c$ and consequently $e(A,A^c)\geq |A|(k-|A|+1)\geq 3(k-2)$. Let $A\subset V$ with $|A|=k-1$. If $|N_1(x)\cap A|\leq k-3$ for every $x\in A$, then $e(A,A^c)\geq 3(k-1)$. Otherwise, let $x\in A$ such that $|N_1(x)\cap A|=k-2$. Denote $N_1(x)\cap A^c=\{y,z\}$. If $\lambda=0$, then each vertex in $N_1(x)\cap A$ has $k-1$ neighbors outside $A$ and thus, $e(A,A^c)\geq 2+(k-2)(k-1)>3k-6$. If $\lambda\geq 1$, then at least $\lambda-1$ of the $\lambda$ common neighbors of $x$ and $y$ are contained in $A$. Therefore, $y$ has at least $\lambda$ neighbors in $A$. A similar statement holds for $z$. Thus, $e(A, N_1(x)\cap A^c) \geq 2\lambda=2(k-b_1-1)$. Also, $e(N_1(x)\cap A, N_2(x))\geq (k-2)b_1$ so $e(A,A^c)\geq (k-2)b_1+2(k-b_1-1)= 3k-6+(k-4)(b_1-1)\geq 3k-6.$
\end{proof}

\begin{theorem}[Tutte \cite{T}] 
A graph $G$ has a perfect matching if and only if $o(G-S)\leq |S|$ for every $S\subset V(G)$.
\end{theorem}

Yu \cite[Theorem 2.2]{Yu} obtained the following characterization of connected graphs that are not $t$-extendable using Tutte's theorem.

\begin{lemma}[Yu \cite{Yu}]\label{nott} 
Let $t\geq 1$ and $G$ be a connected graph containing a perfect matching. The graph $G$ is not $t$-extendable if and only if it contains a subset $S$ of vertices such that the subgraph induced by $S$ contains $t$ independent edges  and $o(G-S)\geq |S|-2t+2$.
\end{lemma}  

The following necessary condition for a bipartite and connected graph not to be $t$-extendable, will be used later in our arguments. 

\begin{lemma}\label{nottbipartite}
Let $G$ be a connected bipartite graph with color classes $X$ and $Y$, where $|X|=|Y|=m$. If $G$ is not $t$-extendable, then $G$ has an independent set $I$ of size at least $m-t+1$, such that $I\not \subset X$ and $I\not \subset Y$. 
\end{lemma}
\begin{proof} Assume that $G$ is not $t$-extendable. By Lemma \ref{nott}, there is a vertex disconnecting set $S$ such that the subgraph induced by $S$ contains at least $t$ independent edges and $o(G-S)\geq |S|-2t+2$. Let $S$ be such a disconnecting set of maximum size. Our key observation is that $G-S$ does not have non-singleton odd components. Indeed, note that any non-singleton odd component of $G-S$ induces a bipartite graph with color classes $A$ and $B$. Since $|A|+|B|$ is odd, we get that $|A|\neq |B|$ and assume that $|A|>|B|$. If $S^\prime=S\cup B$, then $S^\prime$ is a vertex disconnecting set with $|S^\prime|>|S|$ and $o(G-S')\geq |S'|-2t+2$, contradicting to the maximality of $|S|$. By a similar argument, $G-S$ contains no even components. Let $I=V(G)\setminus S$. Then $I$ is an independent set of size at least $m-t+1$ since $|I|+|S|=2m$ and $|I|\geq |S|-2t+2$. Assume that $I\subset X$. Then $S$ induces a bipartite graph with one partite set of size at most $t-1$. This makes it impossible for the subgraph induced by $S$ to contain $t$ independent edges. 
\end{proof}

Note that the study of such independent sets in regular bipartite graphs has been done by other authors in different contexts (see \cite{DeWSV} for example).

\begin{lemma}[Lemma 6 \cite{CL}]\label{1connected}
If $G$ is a distance-regular graph of diameter $D\geq 3$, then for any $x\in V(G)$, the subgraph induced by the vertices at distance $2$ or more from $x$, is connected.
\end{lemma}
\begin{proof} As $G$ has diameter $D\geq 3$, then there are 4 vertices, which induce a $P_4$. It is known that $P_4$ has spectrum $\{\frac{1+\sqrt{5}}{2}, \frac{-1+\sqrt{5}}{2}, \frac{1-\sqrt{5}}{2}, \frac{-1-\sqrt{5}}{2}\}$. By eigenvalue interlacing \cite[Corollary 2.5.2]{BH2}, $\lambda_2\geq \frac{-1+\sqrt{5}}{2}>0$. Cioab\u{a} and Koolen \cite[Theorem 3]{CK} proved that if the entry $u_{i-1}$ of the standard sequence $(u_0,u_1,\dots,u_D)$ corresponding to $\lambda_2$, is positive, then for all $x \in V(G)$, $\Gamma_{\geq i}(x)$ is connected, where $\Gamma_{\geq i}(x)$ is the graph induced by the vertex set at distance at least $i$ to vertex $x$. 
As $u_1=\lambda_2/k>0$, the conclusion follows. \end{proof}

\begin{lemma}[Brouwer and Haemers \cite{BH}]\label{shadowlemma}
Let $G$ be a distance-regular graph. Let $T$ be a disconnecting set of edges of $G$, and let $A$ be the vertex set of a component of $G-T$. Fix a vertex $a \in A$ and let $t_i$ be the number of edges in $T$ that join $\Gamma_{i-1}(a)$ and $\Gamma_i(a)$. Then $|A\cap \Gamma_i(a)|\geq (1-\sum_{j=1}^i\frac{t_j}{c_jk_j})k_i$ and 
$$|A|\geq v-\sum_i\frac{t_i}{c_ik_i}(k_i+\dots +k_D).
$$
If $T$ is a disconnecting set of edges none of which is incident with $a$, then
$$ |A|>v\left ( 1-\frac{|T|}{\mu k_2}\right ).
$$
\end{lemma}

\begin{lemma}\label{Independentset}
Let $G$ be a distance-regular graph with $\lambda \geq 1$. If $A$ is an independent set of $G$, then 
$|N(A)|\geq 2|A|$.
\end{lemma}
\begin{proof}
For any $x\in N(A)$, $N(x)\cap A$ is an independent set in the subgraph $\Gamma_1(x)$. As $\Gamma_1(x)$ is $\lambda$-regular graph with $k$ vertices, its independence number is at most $k/2$. Thus,  $|N(x)\cap A|\leq k/2$. Therefore, $|A|k=e(A,N(A))=\sum_{x\in N(A)}|N(x)\cap A|\leq |N(A)|k/2$ which implies that $|N(A)|\geq 2|A|$.
\end{proof}

\begin{lemma}\label{Independentset2}
Let $G$ be a distance-regular graph with valency $k\geq 3$, $\lambda \geq 1$ and $\mu\leq k/2$. If $A$ is an independent set of $G$, then $|N(A)|\geq k+|A|-1$.
\end{lemma}
\begin{proof}
Let $a=|A|$. The case $a=1$ is trivial. If $a\geq k-1$, Lemma \ref{Independentset} implies that $|N(A)|\geq 2a \geq a+k-1$. Assume that $2\leq a \leq k-2$. If there are two vertices $x, y \in A$, such that $N(x)\cap N(y) =\emptyset$, then $|N(A)|\geq |N(x)\cup N(y)| \geq 2k \geq k+a-1$. Assume that $N(x)\cap N(y) \neq \emptyset$ for any $x,y\in A$. Since $A$ is an independent set, $|N(x) \cap N(y)|=\mu$ for any $x\neq y\in A$. For $z\in N(A)$, let $d_z=|A\cap N(z)|$ and $\bar{d}=\frac{\sum_{z\in N(A)} d_z}{|N(A)|}$. Counting the edges between $A$ and $N(A)$, we have $ak=|N(A)|\bar{d}$. Counting the $3$-subsets of the form $\{x,y,z\}$ such that $x\neq y \in A, z\in N(A), x\sim z, y \sim z$ and then using Jensen's inequality for the function $f(t)={t\choose 2}$, we get that
${a\choose 2}\mu=\sum_{z\in N(A)} {d_z \choose 2}\geq 
|N(A)|{\bar{d} \choose 2}$. Combining these facts, we obtain that $(a-1)\mu\geq k\left(\frac{ka}{|N(A)|}-1\right)$ which implies that $|N(A)|\geq \frac{k^2a}{k+a\mu-\mu}$. As $\mu \leq k/2$, we have $|N(A)| \geq \frac{k^2a}{k+(a-1)k/2}=\frac{2ka}{a+1}=k+a-1+\frac{(a-1)(k-a-1)}{a+1}\geq k+a-1$. 
\end{proof}

A distance-regular graph with intersection array $\{ k, \mu, 1; 1, \mu, k\}$ is called a {\em Taylor graph}. The following lemma due to Brouwer and Koolen (see \cite[Lemma 3.14]{BK2} and also \cite[Proposition 5]{KP} for a generalization) gives a sufficient condition for a distance-regular graph to be a Taylor graph. 
\begin{lemma}[Brouwer and Koolen \cite{BK2}]\label{Taylorgraphbigmu} Let $G$ be a non-bipartite distance-regular graph with $D\geq 3$. If $k< 2\mu$, then $G$ is a Taylor graph.
\end{lemma}

\subsection{Lower bounds for the extendability of distance-regular graphs}\label{lowerbnds}

In this subsection, we give some sufficient conditions, in terms of $k$, $\lambda$ and $\mu$, for a distance-regular graph to be $t$-extendable, where $t\geq 1$.

\begin{theorem}[Chen \cite{C}]\label{Chen}
Let $t\geq 1$ and $n\geq 2$ be two integers. If $G$ is a $(2t+n-2)$-connected $K_{1,n}$-free graph of even order, then $G$ is $t$-extendable. 
\end{theorem}

\begin{corollary}\label{k/4}
If $G$ is a distance-regular graph with even order and $\lambda\geq 1$, then $G$ is $\lfloor \frac{1}{2}\lceil \frac{k+2}{2} \rceil \rfloor$-extendable. 
\end{corollary}
\begin{proof}
The graph $G$ is $K_{1,\lfloor k/2 \rfloor+1}$-free because $\lambda\geq 1$. Let $t=\lfloor \frac{1}{2}\lceil \frac{k+2}{2} \rceil \rfloor$ and $n=\lfloor k/2 \rfloor+1$. Then $k\geq 2t+n-2$. The result follows from Theorem \ref{kconnected} and Theorem \ref{Chen}.
\end{proof}

We improve the previous result when $\mu=1$.

\begin{theorem}
If $G$ is a distance-regular graph with even order, $\lambda\geq 1$ and $\mu=1$, then $G$ is $\lfloor \frac{k+1-\frac{k}{\lambda+1}}{2}\rfloor $-extendable. 
\end{theorem}
\begin{proof}
The condition $\mu=1$ implies that $\Gamma_1(x)$ is a disjoint union of cliques on $\lambda+1$ vertices, for any vertex $x$ of $G$. Hence, $\lambda+1$ divides $k$ and $G$ is $K_{1,\frac{k}{\lambda+1}+1}$-free. Let $t=\lfloor \frac{k+1-\frac{k}{\lambda+1}}{2}\rfloor$ and $n=\frac{k}{\lambda+1}+1$. Then $2t+n-2\leq k$. The conclusion follows from Theorem \ref{kconnected} and Theorem \ref{Chen}.
\end{proof}

The following theorem is an improvement of Corollary \ref{k/4} when $3\leq \mu \leq k/2$.

\begin{theorem}\label{drgn}
Let $G$ be a distance-regular graph with even order, and $D\geq 3$. If $\lambda\geq 1$ and $3\leq \mu \leq k/2$, then $G$ is $t$-extendable, where $t=\lceil \frac{(k-3)(k-1)}{3k-6} \rceil$.
\end{theorem}
\begin{proof}Note that $3\leq \mu \leq k/2$ implies that $k\geq 6$. If $G$ is not $t$-extendable, by Lemma \ref{nott}, there exists a disconnecting $S$ with $s$ vertices such that $o(G-S)\geq s-2t+2$ (and in addition, the subgraph induced by $S$ contains $t$ independent edges). Let $S$ be a disconnecting set with minimum cardinality such that $o(G-S)\geq s-2t+2$. Note that such $S$ may not contain $t$ independent edges. Let $O_1, O_2, \dots, O_r$ be all the odd components of $G-S$, with $r\geq s-2t+2$. Let $a\geq 0$ denote the number of singleton components among $O_1,\dots,O_r$.

We claim that $e(A,S)\geq 3k-6$ for any non-singleton odd component $A$ of $G-S$. 

Let $A$ be a non-singleton odd component of $G-S$ and $B=(A\cup N(A))^c$. If $|A| \leq k-1$, the claim follows from Lemma \ref{outedgeofsmallsets}. Assume that $|A|\geq k$. Let $S^\prime:=\{ s \in N(A) \mid N(s) \subseteq A\cup N(A) \}$. Then $|S^\prime|\leq 1$. Otherwise, assume that  $x\neq y \in S^\prime$. Define $S_0=S\setminus \{x,y\}$ and $A_0=A\cup \{x,y\}$. Then $S_0$ is a disconnecting set with $o(G-S^\prime)=o(G-S)\geq |S|-2t+2>|S^\prime|-2t+2$, contradicting the minimality of $|S|$.  

If we let $A':= \{ a \in A \mid d(a,b) = 2 \text{\ for some \ } b\in B\}$, then $e(A, S) \geq \mu |A^\prime|$. If $|A^\prime| \geq k-2$, we get $e(A,S)\geq \mu|A'|\geq 3(k-2)$ and we are done. Otherwise, if $|A^\prime|<k-2$, then the set $A^\prime \cup S^\prime$ is a disconnecting set with less than $k-1$ vertices, contradicting Theorem \ref{kconnected}. This finish our proof of the claim. 

Counting the number of edges between $S$ and $O_1\cup \dots \cup O_r$, we obtain the following
\begin{equation}
ks\geq e(S,O_1\cup \dots \cup O_r)\geq ak+(r-a)(3k-6)\geq ak+(s-2t+2-a)(3k-6).
\end{equation}
This inequality is equivalent to 
\begin{equation}\label{ineqsa}
t\geq \frac{(k-3)(s-a)+3k-6}{3k-6}
\end{equation}
and since $s-a\geq k-1$ (Lemma \ref{Independentset2}), we obtain that
\begin{equation}
t\geq \frac{(k-3)(k-1)}{3k-6}+1.
\end{equation}
This is a contradiction with $t=\lceil \frac{(k-3)(k-1)}{3k-6}\rceil$.  
\end{proof}
A straightforward calculation shows that $\lceil \frac{(k-3)(k-1)}{3k-6}\rceil= \lfloor \frac {k}{3} \rfloor$ for $k\geq 4$. 
\begin{theorem}\label{Taylorgraph}
Let $G$ be a non-bipartite distance-regular graph with $D\geq 3$ and $\mu > k/2$. Then $\lambda\geq 1$ and $G$ is $t$-extendable, where $t=\lfloor k/3 \rfloor$. 
\end{theorem}   
\begin{proof}
Lemma \ref{Taylorgraphbigmu} implies that $G$ is a Taylor graph with intersection array $\{k, \mu, 1; 1, \mu, k\}$. If $\lambda=0$, then $\mu=k-1$ and $G$ is obtained by deleting a perfect matching from $K_{(k+1)\times(k+1)}$ (see \cite[Corollary 1.5.4]{BCN}) which is a bipartite graph, contradiction.
 
Thus $\lambda\geq 1$. It is known that for any $x\in V(G)$, $\Gamma_1(x)$ is a strongly regular graph with parameters $\left(k,\lambda, \frac{3\lambda-k-1}{2},\frac{\lambda}{2}\right)$ (see \cite[Section 1.5]{BCN}). If $\frac{3\lambda-k-1}{2} \geq 1$, then Lemma \ref{Independentset} implies that $\alpha(\Gamma_1(x))\leq k/3$. If $G$ is not $t$-extendable, then there is a vertex disconnecting set $S$ containing $t$ independent edges, such that  $G-S$ has at least $s-2t+2\geq k-2t+2\geq 3$ odd components. Picking one vertex from each odd component yields an independent set $I$ in $G$. If two vertices of this independent set were at distance $3$, then the neighborhood of these two vertices will be formed by the remaining $2k$ vertices of the graph and therefore, $G-S$ would have only two odd components, contradiction. Thus, any two vertices of this independent set are at distance $2$ to each other. Pick a vertex $x$ in this independent set. Any subset of $k-2t+1$ vertices of $I\setminus \{x\}$ will be an independent set in $\Gamma_1(y)$, where $y$ is the antipodal vertex to $x$. Thus, $k-2t+1\leq k/3$, contradiction with $t=\lfloor k/3\rfloor$. If $\frac{3\lambda-k-1}{2}=0$, then $\Gamma_1(x)$ has parameters $(3\lambda-1,\lambda, 0, \lambda/2)$. If $\lambda=2$, $\Gamma_1(x)$ is $C_5$ which implies that $k=5$ and  $\mu=2$, contradiction with $k/2<\mu$. If $\lambda\geq 4$, then $\Gamma_1(x)$ must have integer eigenvalues implying that $x^2+\frac{\lambda}{2}x-\frac{\lambda}{2}=0$ has integer roots.
However, $(\lambda/2)^2+2\lambda$ is not a perfect square, contradiction. 
\end{proof}

In the end of this subsection, we will show that bipartite distance-regular graphs have high extendability. 

\begin{theorem}\label{bipartitek/2}
If $G$ is a bipartite distance-regular graph with valency $k$, then $G$ is $t$-extendable, where $t=\lfloor \frac{k+1}{2} \rfloor$.  
\end{theorem}
\begin{proof}
Let $X$ and $Y$ be the color classes of $G$, where $|X|=|Y|=m$. Assume that $G$ is not $t$-extendable. By Lemma \ref{nottbipartite}, $G$ has an independent set $I$ of size at least $m-t+1$, such that $I\not \subset X$ and $I\not \subset Y$. Let  $A=I\cap X$, $B=I\cap Y$, $C=X\setminus A$, $D=Y\setminus B$. If $|A|=a$, then $|B|\geq m-a-t+1$, $|C|=m-a$ and $|D|\leq a+t-1$. As there are $ak$ edges between $A$ and $D$, and $(a+t-1)k\geq |D|k=e(D, X)=e(A,D)+e(C,D)$, there are at most $(t-1)k$ edges between $C$ and $D$. This implies that $G$ has an edge cut of size at most $(t-1)k$, which disconnects $G$ into two vertex sets $B\cup C$ and $A\cup D$. Without loss of generality, assume that $|A\cup D|\leq m$. By the second part of Lemma \ref{shadowlemma}, we have 
$$
|A\cup D|>v\left (1-\frac{e(A\cup D, B\cup C)}{\mu k_2} \right)\geq 2m\left (1-\frac{(t-1)k}{(k-1)k}\right )\geq 2m(1-1/2)=m,
$$
contradiction with $|A\cup D|\leq m$. 
\end{proof}

\subsection{The $2$-extendability of distance-regular graphs of valency $k\geq 3$}\label{extend2}

Lou and Zhu \cite{LZ} proved that any strongly regular graph of even order is $2$-extendable with the exception of the complete tripartite graph $K_{2,2,2}$ and the Petersen graph. Cioab\u{a} and Li \cite{CL} showed that any strongly regular graph of even order and valency $k\geq 5$ is $3$-extendable with the exception of the complete $4$-partite graph $K_{2,2,2,2}$, the complement of the Petersen graph and the Shrikhande graph (see \cite[page 123]{BH} for a description of this graph).

In this subsection, we prove that any distance-regular graph of diameter $D\geq 3$ is $2$-extendable. By Corollary \ref{k/4}, any distance-regular graph with $\lambda\geq 1$ and $k\geq 5$ is $2$-extendable. Note also that any distance-regular graph of even order having valency $k\leq4$ and diameter $D\geq 3$ must have  $\lambda = 0$ (see \cite{B,BK}). Theorem \ref{bipartitek/2} implies that any bipartite distance-regular graph of valency $k\geq 3$ is $2$-extendable. Thus, we only need to settle the case of non-bipartite distance-regular graphs with $\lambda=0$. We will need the following lemma.
\begin{lemma}\label{v21}
If $G$ is a non-bipartite distance-regular graph with valency $k\geq 5$ and $\lambda=0$, then $\alpha(G)<v/2-1$.
\end{lemma}
\begin{proof}
If $g$ is the odd girth of $G$, then $v>2g$ and Corollary \ref{independencenumber} implies that $\alpha(G) \leq \frac{v}{2}(1-\frac{1}{g})<v/2-1$. 
\end{proof}

\begin{theorem}\label{maxindependent}
If $G$ is a non-bipartite distance-regular graph with  even order, $D\geq 3$, valency $k\geq 3$ and $\lambda=0$, then $G$ is $2$-extendable.
\end{theorem}
\begin{proof}
We prove this result by contradiction and the outline of our proof is the following. We assume that $G$ is not $2$-extendable. Lemma \ref{nott} implies that there is a vertex disconnecting set $S$, such that the graph induced by $S$ contains at least $2$ independent edges and $o(G-S)\geq |S|-2$. Without loss of generality, we may assume that $S$ is such a disconnecting set with the maximum size. We then prove that $G-S$ does not have non-singleton components which implies that $V(G)-S$ is an independent set of size at least $v/2-1$, contradiction to Lemma \ref{v21}. 

Assume $k\geq 5$ first.

Note that any odd non-singleton component of $G-S$ is not bipartite. Otherwise, assume there is a bipartite odd component of $G-S$ with color classes $X$ and $Y$ such that $|X|>|Y|$. Let $S^\prime=S\cup Y$. Then $|S^\prime|>|S|$ and $o(G-S')\geq |S'|-2$, contradiction with $|S|$ being maximum. Also, $G-S$ has no even components. Otherwise, we can add one vertex of one such even component to $S$ and creating a larger disconnecting set and an extra odd component, contradicting again the maximality of $|S|$. It is easy to see that $G-S$ does not have any components with $3$ vertices, because $G$ is triangle free and any component with $3$ vertices must be a path, hence bipartite.

Assume that $A$ is an odd non-singleton component of $G-S$. If we can show that $e(A,S)\geq 3k-3$, then using $o(G-S)\geq |S|-2$ and $e(X,S)\geq k$ for any component $X$ of $G-S$ (from Theorem \ref{kedge}), we obtain the following contradiction by counting the edges between $S$ and $S^c$:
\begin{equation}\label{edgecount}
k|S|-4\geq e(S,S^c)\geq 3k-3+k(|S|-3)=k|S|-3,
\end{equation}
finishing our proof.

We now prove $e(A,S)\geq 3k-3$ whenever $A$ is a non-singleton odd component of $G-S$.

If $5\leq |A| \leq 2k-3$, then as $A$ has no triangle, Tur\'{a}n\rq{s} theorem implies that $e(A)\leq \frac{|A|^2-1}{4}$, where $e(A)$ denotes the number of edges with both endpoints in $A$. Thus, $e(A,S)\geq k|A|-2e(A)\geq k|A|-\frac{|A|^2-1}{2}\geq 3k-4$. The last equality is attained when $A$ induces a bipartite graph $K_{k-1,k-2}$. This is impossible as the graph induced by $A$ is not bipartite. Hence, $e(A,S)\geq 3k-3$.

Let $A$ be an odd component of $G-S$ such that $|A|\geq 2k-1$. If every vertex of $A$ sends at least one edge to $S$, then we have two subcases: $\mu\geq 2$ and $\mu=1$.

If $\mu\geq 2$, then we can define $S^\prime:=\{s \in N(A) \mid N(s) \subseteq A\cup N(A) \}$. If $|S^\prime| \geq 3$, then $e(A,S)+2e(S) \geq 3k+1$. This is because $e(A,S)+2e(S)=\sum_{x\in S} |N(x) \cap (A \cup S)|$. As the graph induced by $S$ contains at least 2 independent edges, the previous sum contains at least $4$ positive terms, and at least $3$ of such terms are equal to $k$. On the other hand, as in \eqref{edgecount}, counting the number of edges between $S$ and $S^c$, we get that $e(A,S)+(|S|-3)k \leq e(S,S^c)=|S|k-2e(S)$. Thus, $e(A,S)+2e(S) \leq 3k$, contradiction. If $|S^\prime| \leq 2$, then let $B=(A\cup N(A))^c$ and $A^\prime= \{a \in A \mid \exists b \in B \text{\ \ such that \  } d(a,b) = 2\}$. If $A^\prime=A$, then $|A^\prime|\geq 2k-1\geq k-2$. If $A\neq A^\prime$, then because $A'\cup S'$ is a disconnecting set, Lemma \ref{kconnected} implies that $|A^\prime \cup S^\prime|\geq k$ and therefore, $|A^\prime| \geq k-2$. As each vertex in $A^\prime$ sends at least $\mu$ edges to $S$ and $\mu\geq 2$, we get that $e(A,S) \geq 2k-1+(k-2)(\mu-1) \geq 3k-3$. 

If $\mu=1$, then the graph induced by $A$ contains no triangles and four-cycles. If $|A| \geq 3k-3$, then $e(A,S) \geq 3k-3$, as every vertex of $A$ sends at least one edge to $S$. If $|A| \leq 3k-4$, then $e(A)\leq \frac{|A|\sqrt{|A|-1}}{2}$ since the graph induced by $A$ contains no triangles and four-cycles (see \cite[Theorem 2.2]{GKL} or \cite[Theorem 4.2]{VW}). Since also $2k-1\leq |A|\leq 3k-4$, we get that $e(A,S)=k|A|-2e(A)\geq |A|(k-\sqrt{|A|-1})\geq (2k-1)(k-\sqrt{3k-5})\geq 3k-3$.

The only case remaining is when $|A|\geq 2k-1$ and $A$ has a vertex $x$ having no neighbors in $S$ (such a vertex is called a deep point in \cite{BK}). Note that $A^c$ always has  a deep point because every vertex in $V(G)\setminus (A\cup S)$ is a deep point of $A^c$. We have two cases: 

\begin{enumerate} 

\item When $k\geq 6$, we will show that $e(A,S)\geq 3k-3$. Otherwise, by Lemma \ref{shadowlemma}, 
\begin{equation}\label{sizeofA2}
|A|>v\left (1-\frac{3k-4}{\mu k_2}\right)=v\left (1-\frac{3k-4}{k(k-1)}\right)\geq v/2. 
\end{equation}
The last inequality is true since $k\geq 6$. As $A^c$ always has  a deep point, by Lemma \ref{shadowlemma} again, we get that $|A^c|>v/2$, contradiction.

\item When $k=5$, we do not have inequality \eqref{sizeofA2} so we need a different proof. If $\mu\geq 3>k/2$, by Theorem \ref{Taylorgraph}, $\lambda\geq 1$, contradiction. So, we must have $1\leq \mu \leq 2$.

We first show that $A$ is the only non-singleton component of $G-S$. Assume that there are at least two non-singleton components in $G-S$. Let $B$ be another non-singleton component of $G-S$. Then $B$ has a deep point, by previous arguments. If $e(A,S) \geq 2k-1$ and $e(B,S)\geq 2k-1$, then $k|S|-4\geq e(S,S^c)\geq 2(2k-1)+(|S|-4)k=k|S|-2$, contradiction. Without loss of generality, assume that $e(A,S) \leq 2k-2$. By Lemma \ref{shadowlemma}, 
$|A|>v\left (1-\frac{2k-2}{\mu k_2}\right)=v\left (1-\frac{2}{k}\right)=\frac{3v}{5}$. On the other hand, Lemma \ref{shadowlemma} also implies that $|A^c|>\frac{3v}{5}$, contradiction. 

Thus, $A$ is the only non-singleton component in $G-S$. Recall that $|A|\geq 2k-1$ and $A$ has  a deep point $x$. If $e(A,S)\leq 3k-5$, by Lemma \ref{shadowlemma},  $|A|>v\left (1-\frac{3k-5}{\mu k_2}\right)=v\left (1-\frac{10}{20}\right)\geq v/2$. Lemma \ref{shadowlemma} also implies that $|A^c| >v/2$, contradiction. If $e(A,S)=3k-4=11$, by counting the edges between $S$ and $S^c$, we know that $S$ contains exactly two independent edges. Also, $o(G-S)=|S|-2$. Let $X$ be the set of singleton components of $G-S$. We have $|X|=|S|-3$. By Theorem \ref{kconnected}, $|S|\geq k+1 =6$ and $|X|\geq 3$. 

Now, we have two subcases: 

\begin{enumerate}[(i)]

\item Assume that $\mu=2$. Let $W=\{a\in A \mid \exists s\in S, a\sim s\}$. Note that $W\subset A$ and $W$ is a disconnecting set of $G$. By Theorem \ref{kconnected}, $|W|\geq 5$ and the only disconnecting sets of 5 vertices are the neighbors of some vertex.  If $|W|=5$, then we have $W=N(x)$ for some vertex $x$. By Lemma \ref{1connected}, the subgraph induced by the vertices at distance $2$ or more from $x$ is connected. In other word, $W$ disconnects $G$ into two components, $x$ and $V\setminus (W\cup \{x\})$. Since $|A^c|>1$, we must have $A^c=V\setminus (W\cup \{x\})$ and $A\setminus W=\{x\}$. Hence, $|A|=6$, contradicting to $|A|$ is odd. So, $|W|\geq 6$. 

We claim that for any $x \in W$, there exists $t\in X$ such that $d(x,t)=2$. Assume otherwise. Then there is $s \in S$ such that $N(s)\subset A\cup S$. Since the graph induced by $S$ contains exactly two independent edges, $s$ has at most one neighbor in $S$ and at least four neighbors in $A$.  If we let $A^\prime=A\cup \{s\}$ and $S^\prime=S\setminus \{s\}$, then $e(A^\prime, S^\prime)\leq 8$. By Lemma \ref{shadowlemma}, $|A^\prime|>v\left (1-\frac{8}{\mu k_2}\right) =v\left ( 1-\frac{8}{20}\right )=\frac{3v}{5}$.  On the other hand, Lemma \ref{shadowlemma} also implies that $|(A^\prime)^c|>\frac{3v}{5}$, contradiction.  

As $\mu=2$, each vertex in $W$ has at least 2 neighbors in $S$ and $e(A,S)\geq 12$, which is also a contradiction.   

\item Assume that $\mu=1$. We will first prove that $a_2 \leq 1$. If for every $s \in S$, $|N(s)\cap X|\leq 2$, by counting the edges between $S$ and $X$, we have $5|X|=e(S,X) \leq 2|S|$. On the other hand, $|X|=|S|-3 \geq \frac{5}{2}|X|-3$, thus $|X| \leq 2$, contradicting to $|X|\geq 3$. Hence, there exists $s\in S$ such that $|N(s)\cap X|\geq 3$. Let $x,y,z \in N(s) \cap X$. As $\mu=1$, $N(x) \cap N(y) =N(y)\cap N(z)=N(x) \cap N(z) =\{s\}$. Let $U=(N(x)\cup N(y) \cup N(z)) \setminus \{s\}$. It is easy to check that $U\subset N_2(s)$, $|U|=12$, $|N_2(s)|=20$, and $\Gamma_2(s)$ is $a_2$-regular. Since there are at most two edges inside $U$, $12a_2 -4\leq  e(U, N_2(s)\setminus U) \leq 8a_2$ and thus $a_2 \leq 1$.  

Note that $\mu=1$ and $a_2 \leq 1$ imply that $b_2\geq 3$. If there exists $r \in S$, such that $N(r) \subset X$, then $d(r,A)\geq 3$. By Lemma \ref{shadowlemma},  $|A^c|>v\left (1-\frac{3k-4}{k_2b_2}\right)\geq v\left (1-\frac{11}{60}\right)=\frac{49v}{60}$. On the other hand, Lemma \ref{shadowlemma} also implies that $|A|>\left (1-\frac{3k-4}{\mu k_2}\right)=v\left (1-\frac{3k-4}{k(k-1)}\right)=\frac{9v}{20}$, contradiction. Thus, for all $r\in S$, we have  $N(r)\not \subset X$. Consider the edges between $X$ and $S$. We have $5|X|=e(X,S)\leq 4|S|$ and therefore, $|X|=|S|-3\geq 5|X|/4-3$. Thus, $|X|\leq 12$, $|S|\leq 15$, $27\geq |A^c|>9v/20$ and $v<60$. Note that there is no distance-regular graph with $v<60$, $k=5$, $\lambda=0$, $\mu=1$ and $a_2\leq 1$, see the table \cite{B} and \cite{Flaass} (where it was shown that there exists no distance-regular graph with intersection array $\{5,4,3;1,1,2\}$).

\end{enumerate}

\end{enumerate}

This finishes the proof of the case $k\geq 5$.

When $k=4$, all the distance-regular graph with even order are bipartite \cite{BK} so we are done by Theorem \ref{bipartitek/2}.

When $k=3$, there are $3$ non-bipartite triangle-free distance-regular graphs with even order (see \cite{BBT} or \cite[Chapter 7]{BCN}): the Coxeter graph (intersection array $\{3,2,2,1;1,1,1,2\}$), the Dodecahedron graph (intersection array $\{3,2,1,1,1;1,1,1,2,3\}$) and the Biggs-Smith graph (intersection array $\{3,2,2,2,1,1,1,;1,1,1,1,1,1,3\}$). We will show that each one of them is $2$-extendable. 

Let $G$ be the Coxeter graph. Then $G$ has $28$ vertices, girth $7$ and independence number $12$ (see \cite{Biggs} for example). If $G$ is not $2$-extendable, there is a disconnecting set  $S$  of maximum size, such that the graph induced by $S$ contains $2$ independent edges and $o(G-S)\geq |S|-2$.  As $|S|\geq 4$, we have $o(G-S)\geq 2$. Assume that $G-S$ contains a non-singleton component $A$. As $|A^c|\geq |S|+1\geq 5$, we have that $3\leq |A| \leq v-|A^c| \leq 23$. If $3\leq |A| \leq 5$, the graph induced by $A$ is bipartite as the girth of $G$ is $7$. As in the case $k\geq 5$, we can construct a larger disconnecting set contradicting the maximality of $S$. If $|A^c|=5$, then we have that $|S|=4$ and there is one singleton component $\{x\}$ in $A^c$. Since $S$ contains two independent edges and $x$ has three neighbors in $S$, we obtain $\lambda\neq 0$, contradiction.


Let $G$ be the Dodecahedron graph. Then $G$ has $20$ vertices, girth $5$ and independence number $8$ (see \cite[pp.116]{GR} for example). If $G$ is not 2-extendable, there is a disconnecting set $S$ of maximum size, such that the graph induced by $S$ contains $2$ independent edges and $o(G-S) \geq |S|-2$. As $|S|\geq 4$, we have $o(G-S)\geq 2$. Assume that $G-S$ contains a non-singleton component $A$. As  $|A^c|\geq |S|+1\geq 5$, we have $3\leq |A| \leq 15$. We will prove that $|A|\not =3, 5,7,9$ and $|A^c|\not= 5,7,9$. By maximality of $|S|$, the graph induced by $A$ is not bipartite. So, $|A|\not =3$.  If $|A|=7$, then the graph induced by $A$ contains at most one cycle. Thus, $e(A)\leq 7$ and $e(A,A^c)= 3|A|-2e(A)\geq 7$.  If $|A|=9$, then the graph induced by $A$ contains at most two cycles. Thus, $e(A) \leq 10$ and $e(A,A^c) = 3|A|-2e(A)\geq 7$. In either case, we will obtain a contradiction by inequality \eqref{edgecount}. Using the same argument, we can show that $|A^c|\not=7, 9$. If $|A^c|=5$, then we have that $|S|=4$ and there is one singleton component $\{x\}$ in $A^c$. Since $S$ contains two independent edges and $x$ has three neighbors in $S$, we obtain $\lambda\neq 0$, contradiction.


Let $G$ be the Biggs-Smith graph. Then $G$ has girth $9$ and 102 vertices.  If $G$ is not 2-extendable, there is a disconnecting set $S$ of maximum size, such that the graph induced by $S$ contains 2 independent edges and $o(G-S)\geq |S|-2$. Assume that $G-S$ contains a non-singleton component $A$. By similar argument as the previous cases, we can assume that $5\leq |A| \leq 97$. When $5\leq |A| \leq 7$, $e(A)=|A|-1$ and $e(A,A^c)= 3|A|-2e(A)=|A|+2 \geq 7$ . When $9\leq |A| \leq 15$, $e(A)\leq |A|$ and $e(A,A^c)=3|A|-2e(A)\geq |A|\geq 9$.  When $17 \leq |A| \leq 51$, $e(A,A^c)\geq \frac{(3-2.56155)|A|(102-|A|)}{102}\geq 6.21134>6$ (see \cite[Corollary 4.8.4]{BH2} or \cite{Mohar}). If $e(A,S)\geq 3k-3=6$, we will obtain a contradiction by inequality \eqref{edgecount}. Using the same argument, we can obtain a contradiction when $5\leq |A^c|\leq 51$. Thus, all the components of $G-S$ are singletons. Therefore, $\alpha(G)\geq o(G-S)\geq \max\{102-|S|,|S|-2\}\geq 50$, contradiction with $\alpha(G)=43$ (see the table on page 6).
\end{proof}

\section{Final Remarks}

Note that some of the bounds in this paper may be improved if one obtains better lower bound for $e(A,A^c)$ with $k\leq |A| \leq v-k$. We make the following conjecture which is still open for strongly regular graphs \cite{CL}.
\begin{conj}
If $G$ is a  distance-regular graph of valency $k$, even order $v$ and diameter $D\geq 3$, then the extendability of $G$ is at least $\lceil k/2\rceil -1$.
\end{conj}

A stronger property than $m$-extendability is the property $E(m,n)$ introduced by Porteous and Aldred \cite{PA}. A connected graph with at least $2(m+n+1)$ vertices is said to be $E(m,n)$ if for every pair of disjoint matchings $M,N$ of $G$ of size $m$ and $n$, respectively, there exists a perfect matching in $F$ such that $M\subseteq F$ and $F\cap N=\emptyset$. It would be interesting to investigate this property for distance-regular graphs and graphs in association schemes. Godsil \cite{Godsil} conjectured that the edge-connectivity of a connected class of an association scheme equals its valency and Brouwer \cite{Brouwer} made the stronger conjecture that the vertex-connectivity equals the valency. Brouwer's conjecture has been proved by Brouwer and Koolen \cite{BK2} for distance-regular graphs, but both Godsil and Brouwer's conjectures are open in the other cases. Godsil's conjecture would imply that any connected class in an association scheme of even order, has a perfect matching. To our knowledge, this is not known at present time.

\section*{Acknowledgements} We thank Bill Martin and the two anonymous referees for many useful comments and suggestions that have greatly improved our initial manuscript.

\end{document}